\documentclass{amsart}
\usepackage{amssymb, latexsym}
\DeclareMathOperator{\im}{Im}

\DeclareMathOperator{\re}{Re}

\newcommand{\eps}{{\varepsilon}}
\newcommand{\alg}{{\mathfrak A}}
\newcommand{\seqalg}{{\widehat \alg}}
\newcommand{\algloc}{{\alg_{loc}}}

\newcommand{\D}{{\mathbb D}}

\newcommand{\C}{{\mathbb C}}
\newcommand{\R}{{\mathbb R}}
\newcommand{\N}{{\mathbb N}}

\newcommand{\sm}{{\mathcal M}}
\newcommand{\interval}{{[-2,2]}}
\newcommand{\V}{{\mathcal{V}}}
\newcommand{\Szego}{{Szeg\H{o} }}

\theoremstyle{plain}
\newtheorem{theorem}{Theorem}
\newtheorem{lemma}[theorem]{Lemma}

\newtheorem{proposition}[theorem]{Proposition}

\numberwithin{equation}{section} \numberwithin{theorem}{section}

\begin{document}

\title[A spectral equivalence for Jacobi matrices]{A spectral equivalence for Jacobi matrices}
\author{E. ~Ryckman}
\address{University of California, Los Angeles CA 90095, USA}
\email{eryckman@math.ucla.edu}

\begin{abstract}
We use the classical results of Baxter and Gollinski-Ibragimov to
prove a new spectral equivalence for Jacobi matrices on $l^2(\N)$.
In particular, we consider the class of Jacobi matrices with
conditionally summable parameter sequences and find necessary and
sufficient conditions on the spectral measure such that
$\sum_{k=n}^\infty b_k$ and $\sum_{k=n}^\infty (a_k^2 - 1)$ lie in
$l^2_1 \cap l^1$ or $l^1_s$ for $s \geq 1$.
\end{abstract}

\maketitle

%%%%%%%%%%%%%%%%%%%%%%%%%%%%%%%%%%%%%%%%%%%%%%%%%%%%%%%%%%%%%%%%%%%%%%%%%%%%%%%%%%%%%%%%%%%
%
%
%                                   Section
%
%
%%%%%%%%%%%%%%%%%%%%%%%%%%%%%%%%%%%%%%%%%%%%%%%%%%%%%%%%%%%%%%%%%%%%%%%%%%%%%%%%%%%%%%%%%%%

\section{Introduction}
Let us begin with some notation.  We study the spectral theory of
Jacobi matrices, that is semi-infinite tridiagonal matrices
$$J = \begin{pmatrix}
b_1 & a_1 & 0 & 0 \\
a_1 & b_2 & a_2 & 0 \vphantom{\ddots}\\
0 & a_2 & b_3 & \ddots\\
0 & 0 & \ddots & \ddots
\end{pmatrix}$$
where $a_n > 0$ and $b_n \in \R$.  In this paper we make the
overarching assumption that the sequences $b_n$ and $a_n^2-1$ are
conditionally summable. We may then define
\begin{align}\label{lambda kappa def}
\begin{split}
\lambda_n &:= -\sum_{k=n+1}^\infty b_k\\
\kappa_n &:= -\sum_{k=n+1}^\infty (a_k^2 - 1)
\end{split}
\end{align}
for $n = 0,1,\dots$.

Let $d\nu$ be the spectral measure for the pair $(J,\delta_1)$,
where $\delta_1 = (1,0,0,\dots)^t$, and assume that $d\nu$ is not
supported on a finite set of points (we call such measures
\emph{nontrivial}). Let
\begin{equation*}
m(z) := \langle \delta_1 , (J-z)^{-1}\delta_1 \rangle = \int
\frac{d\nu(x)}{x-z}
\end{equation*}
be the associated $m$-function, defined for $z \in \C \backslash
\text{supp}(\nu)$.

Recall that
$$\{\beta_n\} \in l^p_{s} \quad\text{if}\quad \|\beta\|_{l^p_{s}}^p := \sum_n
|n|^s|\beta_n|^p < \infty.$$ Throughout, let $\seqalg$ denote
either of the algebras $l^2_1 \cap l^1$ or $l^1_s$ where $s \geq
1$, and $\alg$ the set of (complex valued) functions on the circle
$\partial \D$ whose Fourier coefficients lie in $\seqalg$. Notice
that every $f \in \alg$ has $l^1$ Fourier coefficients so is
continuous.  If $f$ is a function on $\interval$, we write $f \in
\alg$ if $f(2 \cos \theta) \in \alg$.  Finally, we will say that
$d\nu \in \V$ if
\begin{enumerate}
\item $J$ has finitely-many eigenvalues and they all lie in $\R
\setminus \interval$

\item $d\nu$ is absolutely continuous on $\interval$ and may be
written there as
$$\bigl( \sqrt{2+x} \bigr)^l \bigl( \sqrt{2-x} \bigr)^r
v_0(x) dx$$ where $l,r \in \{\pm 1\}$ and $\log v_0 \in \alg$.
\end{enumerate}

Our main result is\footnote{A similar result is proved in
\cite{Ryckman}, but with $\seqalg$ replaced by $l^2_1$.  While the
techniques of that paper extend to handle the case discussed here,
they are quite lengthy and involved.  Our aim is to provide a
proof of this simpler result that is both general and short.}:
\begin{theorem}\label{lemma main thm alg case}
Let $J$ be a Jacobi matrix.  The following are equivalent:
\begin{enumerate}
\item The sequences associated to $J$ by \eqref{lambda kappa def}
obey $\lambda,\kappa \in \seqalg$

\item $d\nu \in \V$.
\end{enumerate}
\end{theorem}

The main ingredient in the proof will be the following versions of
the Strong \Szego Theorem and Baxter's Theorem\footnote{The
version of the Strong \Szego Theorem we use is due to
\cite{Golinskii Ibragimov} and \cite{Ibragimov}.  The version of
Baxter's Theorem is due to \cite{Simon}.  For relevant definitions
see, for instance, \cite{Simon}.}.

\begin{theorem}[Golinskii-Ibragimov]\label{lemma Golinskii-Ibragimov}
Let $d\mu$ be a nontrivial probability measure on $\partial \D$
with Verblunsky parameters $\{\alpha_n\} \subseteq \D$. The
following are equivalent:
\begin{enumerate}
\item $\alpha \in l^2_1$

\item $d\mu = w \frac{d\theta}{2\pi}$ and $(\log w)^\wedge \in
l^2_1$.
\end{enumerate}
\end{theorem}

\begin{theorem}[Baxter]\label{lemma baxter}
Let $d\mu(\theta)$ be a nontrivial probability measure on
$\partial \D$ with Verblunsky parameters $\{\alpha_n\}$, and let
$s \geq 0$. The following are equivalent:
\begin{enumerate}
  \item $\alpha \in l^1_s$

  \item $d\mu = w \frac{d\theta}{2\pi}$ and $(\log w)^\wedge \in
l^1_s$.
\end{enumerate}
\end{theorem}

In Section 2 we develop relations between the Jacobi and
Verblunsky parameters, in Section 3 we discuss the relationship
between measures and $m$-functions, in Section 4 we prove some
results about adding and removing eigenvalues, and in Section 5 we
prove Theorem \ref{lemma main thm alg case}. To motivate the
results of Sections 2 and 4 we outline the proof now. Let
$J^{(N)}$ be the Jacobi matrix obtained from $J$ by removing the
first $N$ rows and columns. To prove the forward direction, we
choose $N$ large enough that $\sigma (J^{(N)}) \subseteq [-2,2]$
so the Verblunsky parameters exist. The results of Section 2 then
allow us to apply Theorems \ref{lemma Golinskii-Ibragimov} and
\ref{lemma baxter} to an operator differing from $J^{(N)}$ in the
first row and column to see that this operator has a spectral
measure with the correct form. We conclude the proof by using the
results of Sections 3 and 4 to show that the conditions on the
spectral measure are unaffected by changing the top row and column
of the operator, or by adding back on the removed rows and
columns.  To prove the reverse implication we essentially run this
argument backward.

It is a pleasure to thank Rowan Killip for his helpful advice.

%%%%%%%%%%%%%%%%%%%%%%%%%%%%%%%%%%%%%%%%%%%%%%%%%%%%%%%%%%%%%%%%%%%%%%%%%%%%%%%%%%%%%%%%%%%
%
%
%                                   Section
%
%
%%%%%%%%%%%%%%%%%%%%%%%%%%%%%%%%%%%%%%%%%%%%%%%%%%%%%%%%%%%%%%%%%%%%%%%%%%%%%%%%%%%%%%%%%%%

\section{The Geronimus Relations}
Given a nontrivial probability measure $d\mu$ on $\partial \D$
that is invariant under conjugation, define a nontrivial
probability measure $d\nu$ on $[-2,2]$ by
$$\int_{-2}^2 g(x) d\nu(x) = \int_0^{2\pi} g(2 \cos \theta) d\mu
(\theta).$$ Similarly, given such a measure $d\nu$, one can define
a measure $d\mu$ by $$\int_0^{2\pi} h(\theta) d\mu(\theta) =
\int_{-2}^2 h(\arccos (x/2)) d\nu(x)$$ when $h(-\theta) =
h(\theta)$. It is clear that $d\mu$ is a nontrivial probability
measure that is invariant under conjugation.

This sets up a one-to-one correspondence between the set of
nontrivial probability measures on $[-2,2]$ and the set of
nontrivial probability measures on $\partial \D$ invariant under
conjugation. We call the map $d\mu \mapsto d\nu$ the Szeg\H{o}
mapping and denote it by $d\nu = Sz(d\mu)$. If the two measures
are absolutely continuous with respect to Lebesgue measure we will
write $d\mu(\theta) = w(\theta) \frac{d\theta}{2\pi}$ and $d\nu(x)
= v(x) dx$.  In this case we have
\begin{equation}\label{weight relation}
\begin{split}
w(\theta) = 2\pi |\sin(\theta)| v(2\cos(\theta))\\
v(x) = \frac{1}{\pi \sqrt{4-x^2}}w(\arccos(x/2)).
\end{split}
\end{equation}

The connection between $\alpha$, and $a,b$ is given by

\begin{theorem}[The Geronimus Relations \cite{Geronimus}]\label{lemma Geronimus Relations}
Let $d\mu$ be a nontrivial probability measure on $\partial \D$
that is invariant under conjugation, and let $d\nu = Sz(d\mu)$.
 Then for all $n \geq 0$
\begin{align}\label{Geronimus Relations}
\begin{split}
a_{n+1}^2 &= (1-\alpha_{2n-1})(1-\alpha_{2n}^2)(1+\alpha_{2n+1})\\
b_{n+1} &= (1-\alpha_{2n-1})\alpha_{2n} -
(1+\alpha_{2n-1})\alpha_{2n-2}.
\end{split}
\end{align}
\end{theorem}

Since $a_n > 0$, there is no ambiguity in which sign to choose for
the square root in \eqref{Geronimus Relations}.  Unless otherwise
noted we take $\alpha_{-1} = -1$.  The value of $\alpha_{-2}$ is
irrelevant since it is multiplied by zero.

From \eqref{Geronimus Relations} we see that decay of the
$\alpha$'s determines decay of the $a$'s and $b$'s. However, given
sequences $a, b$ it is difficult to determine whether the
corresponding $\alpha$ sequence even exists\footnote{The existence
of $\alpha$ is equivalent to $\sigma (J) \subseteq [-2,2]$.  See,
for instance, \cite{Damanik Killip}.}, and then whether decay of
$a, b$ is passed to $\alpha$. The rest of this section is devoted
to resolving these problems. We begin with the simple

\begin{lemma}\label{lemma seq splitting} Let $p = 1,2$, $s \geq
1$, $\beta, \gamma \in l^p_{s}$, and define a sequence $\eta_n :=
\sum_{k=n}^\infty \beta_k \gamma_k$.  Then\footnote{We do not
expect such a result for $l^p_s$ if $0 \leq s < 1$, as can be seen
by considering $\beta_n = \gamma_n = \delta_N (n)$.} $\eta \in
l^p_{s}$ and $\|\eta\|_{l^p_{s}} \leq
\|\beta\|_{l^p_{s}}\|\gamma\|_{l^p_{s}}$.
\end{lemma}

\begin{proof}
First consider the $l^2_s$ case.  By Cauchy-Schwarz and the
definition of $\|\cdot\|_{l^2_s}$ we have
\begin{align*}
\|\eta\|^2 &= \sum_{n=1}^\infty
n^{s}\Bigl|\sum_{k=n}^\infty \beta_k \gamma_k \Bigr|^2 \leq
\sum_{n=1}^\infty n^{s} \Bigl(\sum_{k=n}^\infty |\beta_k
\gamma_k|\Bigr)^2 \\
&\leq \sum_{n=1}^\infty n^{s} \Bigl(\sum_{k=n}^\infty |\beta_k|^2
\Bigr) \Bigl(\sum_{k=n}^\infty |\gamma_k|^2 \Bigr) =
\sum_{n=1}^\infty \Bigl(\sum_{k=n}^\infty n^{s} |\beta_k|^2
\Bigr) \Bigl(\sum_{k=n}^\infty |\gamma_k|^2 \Bigr) \\
&\leq \|\beta\|^2 \sum_{n=1}^\infty \sum_{k=n}^\infty |\gamma_k|^2
= \|\beta\|^2 \sum_{n=1}^\infty n |\gamma_n|^2 \leq \|\beta\|^2
\cdot \|\gamma\|^2.
\end{align*}
For $l^1_s$ we replace the use of Cauchy-Schwarz by $$\sum
|\beta_k \gamma_k| \leq \Bigl(\sum |\beta_k|\Bigr)\Bigl(\sum
|\gamma_k|\Bigr)$$ then argue as above.
\end{proof}

We can use Lemma \ref{lemma seq splitting} to show that decay of
the $\alpha$'s is inherited by $\lambda$ and $\kappa$. Given
numbers $x$ and $y$, we will write $x \lesssim y$ if there exists
a $c > 0$ so that $x \leq cy$.

\begin{lemma}\label{lemma alpha implies parameters}
Let $\{\alpha_n\} \subseteq [-1,1]$ and $\alpha \in l^p_s$ for $p
= 1,2$ and $s \geq 1$.  Define
\begin{gather*}
K(\alpha)_n = \sum_{k=n}^\infty \alpha_{2k}^2 +
\alpha_{2k-1}\alpha_{2k+1} - \alpha_{2k}^2(\alpha_{2k-1} -
\alpha_{2k+1}) - \alpha_{2k}^2\alpha_{2k-1}\alpha_{2k+1}\\
 L(\alpha)_n = \sum_{k=n}^\infty
\alpha_{2k-1}(\alpha_{2k} + \alpha_{2k-2}).
\end{gather*}
Then $L(\alpha), K(\alpha) \in l^p_s$ with
\begin{equation*}
\|L(\alpha)\|_{l^p_s} + \|K(\alpha)\|_{l^p_s} \leq C
\|\alpha\|_{l^p_s}^2
\end{equation*}
for some $C > 0$. If $\alpha \in l^2_1 \cap l^1$, then $\lambda,
\kappa \in l^2_1 \cap l^1$ as well.
\end{lemma}

By expanding the right-hand side of \eqref{Geronimus Relations}
one obtains
\begin{equation}\label{expanded lambda kappa}
\begin{split}
\kappa_n = \alpha_{2n-1} + K(\alpha)_n\\
\lambda_n = \alpha_{2n-2} + L(\alpha)_n.
\end{split}
\end{equation}
So the above result may be interpreted as saying that when $\alpha
\in l^p_s$ or $l^1 \cap l^2_s$ then so are $\lambda$ and $\kappa$,
with norms depending on that of $\alpha$.

\begin{proof}
To see the $l^p_s$ statement holds, we use the bound $(|a|+|b|)^2
\leq 2(|a|^2 + |b|^2)$, the hypothesis that $|\alpha_n| \leq 1$,
and repeated applications of Lemma \ref{lemma seq splitting}.

To prove $K \in l^1$ if $\alpha \in l^2_1 \cap l^1$ we use
\begin{align*}
\|K\|_{l^1} &= \sum_{n=0}^\infty \Bigl| \sum_{k=n}^\infty
\alpha_{2k}^2 + \alpha_{2k-1}\alpha_{2k+1} -
\alpha_{2k}^2(\alpha_{2k-1} -
\alpha_{2k+1}) - \alpha_{2k}^2\alpha_{2k-1}\alpha_{2k+1} \Bigr| \\
&\leq \sum_{n=0}^\infty \sum_{k=n}^\infty |\alpha_{2k}^2 +
\alpha_{2k-1}\alpha_{2k+1} - \alpha_{2k}^2(\alpha_{2k-1} -
\alpha_{2k+1}) - \alpha_{2k}^2\alpha_{2k-1}\alpha_{2k+1}| \\
&\leq \sum_{k=1}^\infty (k-1) |\alpha_{2k}^2 +
\alpha_{2k-1}\alpha_{2k+1} - \alpha_{2k}^2(\alpha_{2k-1} -
\alpha_{2k+1}) - \alpha_{2k}^2\alpha_{2k-1}\alpha_{2k+1}| \\
&\lesssim \|\alpha\|_{l^2_s}^2.
\end{align*}
The proof that $L \in l^1$ is similar.
\end{proof}

Given two sequences $\lambda, \kappa$, we now investigate when
there exists a sequence $\alpha$ solving \eqref{expanded lambda
kappa}.  The main step is the following technical bound.

\begin{lemma}\label{lemma Lipschitz condition}
Let $p = 1, 2$, $s \geq 1$, and $\lambda,\kappa \in l^p_s$ be
given. Define a map $F$ by
$$F(\beta)_{2n-1} = \lambda_n + L(\beta)_n
\quad\text{and}\quad F(\beta)_{2n} = \kappa_n + K(\beta)_n.$$ Then
$F : l^p_s \rightarrow l^p_s$, and for any $\beta,\gamma \in
l^p_s$ we have
$$\|F(\beta)-F(\gamma)\|_{l^p_s} \leq C' (\|\beta\|_{l^p_s} +
\|\gamma\|_{l^p_s})^{1/2}\|\beta - \gamma\|_{l^p_s}^2$$ for some
$C'>0$.
\end{lemma}

\begin{proof}
By Lemma \ref{lemma alpha implies parameters}, the range of $F$ is
as stated. Now let $\beta,\gamma \in l^p_s$. We'll bound the sum
over odd values of $n$, the proof for even values of $n$ is
analogous. For $p = 2$ we have
\begin{align*}
\|F(\beta)-F(\gamma)\|^2 &= \|L(\beta) - L(\gamma)\|^2 \\
&= \sum_{n=1}^\infty (2n-1)^{s} \Biggl| \sum_{k=2n-1}^\infty
\beta_{2k-1}(\beta_{2k}
+ \beta_{2k-2}) - \gamma_{2k-1}(\gamma_{2k} + \gamma_{2k-2}) \Biggr|^2 \\
&\lesssim \sum_{n=1}^\infty (2n-1)^{s} \Biggl( \Bigl(
\sum_{k=2n-1}^\infty |\beta_{2k} + \beta_{2k-2}|\cdot|\beta_{2k-1}
- \gamma_{2k-1}|
\Bigr)^2 \\
&\quad+ \Bigl(\sum_{k=2n-1}^\infty
|\gamma_{2k-1}|\cdot|(\beta_{2k} + \beta_{2k-2}) - (\gamma_{2k} +
\gamma_{2k-2})| \Bigr)^2 \Biggr)
\end{align*}
where the inequality was obtained by adding and subtracting the
term $$\gamma_{2k-1}(\beta_{2k} + \beta_{2k-2}).$$  By Lemma
\ref{lemma seq splitting} we can replace the sums in $k$ by norms
to bound
\begin{align*}
\|F(\beta)-F(\gamma)\|^2 & \lesssim \Bigl( \|\beta_{2k} +
\beta_{2k-2}\|^2 \cdot \|\beta_{2k-1} - \gamma_{2k-1}\|^2\\
 &\quad +\|\gamma_{2k-1}\|^2 \cdot \|(\beta_{2k} + \beta_{2k-2}) -
(\gamma_{2k} + \gamma_{2k-2})\|^2 \Bigr)\\
&\lesssim \Bigl(\|\beta\|^2 + \|\gamma\|^2 \Bigr) \|\beta -
\gamma\|^2
\end{align*}
as claimed.  The proof for $p = 1$ is similar and simpler.
\end{proof}

\begin{proposition}\label{lemma solves Geronimus}
Given $\lambda,\kappa \in \seqalg$ with small enough norms, there
exists a sequence $\alpha \in \seqalg$ solving \eqref{expanded
lambda kappa}.
\end{proposition}

\begin{proof}
Let $\|\cdot\|$ denote the norm on $\seqalg$, and let $C$ be the
universal constant arising in Lemma \ref{lemma alpha implies
parameters}: $\|L(\beta)\|+\|K(\beta)\| \leq C \|\beta\|^2$. Let
$0 < \eps < \tfrac{1}{2C}$ to be chosen momentarily, and suppose
that $\|\lambda\|,\|\kappa\| \leq \eps (1/2-C\eps)$. Then by
considering the even and odd terms separately we see
$\|F(\beta)_{\text{odd}}\| \leq \|\lambda\| + \|L(\beta)\| \leq
\eps/2$ if $\|\beta\| \leq \eps$, and similarly
$\|F(\beta)_{\text{even}}\| \leq \eps/2$. So $F$ maps the
$\eps$-ball in $\seqalg$ back to itself.  By Lemma \ref{lemma
Lipschitz condition}, it is Lipschitz on the $\eps$-ball with
Lipschitz constant $\sqrt{2} C' \eps$, where $C'$ is the universal
constant arising in Lemma \ref{lemma Lipschitz condition}. So if
$\eps$ is small enough, the Banach Fixed Point Theorem provides a
unique fixed point $\alpha$ of $F$ with $\|\alpha\| < \eps$. From
the definition of $F$ we see this fixed point solves
\eqref{expanded lambda kappa} with the prescribed $\lambda$ and
$\kappa$.
\end{proof}

%%%%%%%%%%%%%%%%%%%%%%%%%%%%%%%%%%%%%%%%%%%%%%%%%%%%%%%%%%%%%%%%%%%%%%%%%%%%%%%%%%%%%%%%%%%
%
%
%                                   Section
%
%
%%%%%%%%%%%%%%%%%%%%%%%%%%%%%%%%%%%%%%%%%%%%%%%%%%%%%%%%%%%%%%%%%%%%%%%%%%%%%%%%%%%%%%%%%%%

\section{Connecting $d\nu$ and $m$}
In the next section we will begin to add and remove eigenvalues of
$J$.  It is more convenient to recast the criterion on $d\nu$ from
Theorem \ref{lemma main thm alg case} in terms of its associated
$m$-function, which we do in this section.

The map $z \mapsto E(z) := z + z^{-1}$ is a conformal mapping of
$\D$ to $\C \cup \{\infty\} \backslash [-2,2]$ sending $0$ to
$\infty$ and $\pm 1$ to $\pm 2$. Define the $M$-function
associated to $d\nu$ as
$$M(z) = -m(E(z)) = -m(z + z^{-1}) = \int \frac{z d\nu(x)}{1 - xz +
z^2}.$$ We have introduced the minus sign so that $M$ is Herglotz;
that is $M$ is analytic on $\D$ and maps $\C_+$ to itself (as $E
\mapsto z$ maps the upper half-plane to the lower half-disc).

The $M$-function encodes all the spectral information of $J$ (see,
for instance, \cite{Killip Simon, Teschl}). The poles of $M$ in
$(-1,1)$ are related to the eigenvalues of $J$ off $[-2,2]$ by the
map $z \mapsto E$.  We can recover the entries of $J$ from the
continued fraction expansion of $M(z)$ near infinity (see
\cite{Stieltjes})
\begin{equation}\label{continued fraction for m}
M(z) = \cfrac{1}{z + z^{-1} - b_1 - \cfrac{a_1^2}{z + z^{-1} - b_2
- \dots}}.
\end{equation}

We will say a Jacobi matrix is resonant at $E = 2$ if
$$\lim_{z \uparrow 1} |M(z)| = \infty,$$ and we will say $J$
is nonresonant at $E = 2$ otherwise.  We define resonance at $E =
-2$ similarly.  If $J$ is resonant at both $E = -2$ and $E = 2$ we
will say $J$ is doubly-resonant.

Throughout what follows, we make frequent use of a theorem of
Wiener and Levy. For convenience we recall it here (a proof can be
found in \cite{Zygmund}):

\begin{theorem}[Wiener-Levy]\label{lemma wiener levy}
Let $\mathfrak{B}$ be a commutative Banach algebra, $x \in
\mathfrak{B}$, and $F$ analytic in a neighborhood of $\sigma (x)$.
Then $F(x) \in \mathfrak{B}$ can be naturally defined so that $F
\mapsto F(x)$ is an algebra homomorphism of the functions analytic
in a neighborhood of $\sigma (x)$ into $\mathfrak{B}$.
\end{theorem}

Recall that if $f \in \alg$ then its spectrum is its range. So
this shows that if $F$ is analytic in a neighborhood of the range
of $f \in \alg$, then $F(f) \in \alg$ too.

We enrich the algebra framework a bit further by allowing
functions that are only locally in the algebra.  Given $\theta_0
\in [0,2\pi)$, we'll say that $f \in \algloc (\theta_0)$ if there
is a smooth bump $\chi$ on $\partial \D$ equalling one in a
neighborhood of $\theta_0$ such that $\chi f \in \alg$.  Given an
open interval $I \subseteq \partial \D$, we will say $f \in
\algloc (I)$ if $f \in \algloc (\theta_0)$ for all $\theta_0 \in
I$.  Notice that the bumps $\chi$ are in $\alg$, so if $f \in
\algloc (\theta_0)$ for all $\theta_0 \in [0,2\pi)$, then by
choosing a partition of unity on $[0,2\pi)$ with small enough
supports we see $f \in \alg$ too.

We now transfer criterion on $d\nu$ to criterion on $M$. We will
write $M \in \sm$ if
\begin{enumerate}
\item For all intervals $I \subseteq \partial \D$ avoiding $z =
\pm 1$ we have $M \in \algloc (I)$ and $\im M \neq 0$ on $I$.

\item For $z_0 \in \{-1,1\}$, there is a $\partial
\D$-neighborhood $I$ of $z_0$ and a smooth bump $\chi$ supported
on $I$ and equalling one near $z_0$, such that either
$$\chi(\theta) M(\theta) = \chi(\theta)\frac{G(\theta)}{\sin
(\theta)}$$ or
$$\chi(\theta)M(\theta) = \chi(\theta)\Bigl( c + \sin(\theta)G(\theta)
\Bigr)$$ where $c \in \R$, $G \in \alg$, and $\im G \neq 0$ on
$I$.
\end{enumerate}

\begin{proposition}\label{lemma dnu M connection}
For any Jacobi matrix $J$, $d\nu \in \V$ if and only if $M \in
\sm$.
\end{proposition}

In particular, to prove Theorem \ref{lemma main thm alg case}, it
suffices to prove that $\lambda, \kappa \in \seqalg$ if and only
if $M \in \sm$. In the proof of Proposition \ref{lemma dnu M
connection} we will make frequent use of the following lemma.

\begin{lemma}\label{lemma hilbert transform}
Let $H$ be the Hilbert transform on $\partial \D$.  Let $f$ be a
smooth function on $\partial \D$, and let $A_f$ represent the
operator $g(\theta) \mapsto f(\theta)g(\theta).$  If $\eta$ is a
measure on $\partial \D$, then $[A_f,H]\eta$ is a smooth function
(where $[A,B] = AB - BA$ is the usual commutator bracket).
\end{lemma}

We will not prove this here.  It is a fairly standard result from
Harmonic Analysis.

\begin{proof}[Proof of Proposition \ref{lemma dnu M connection}]
Recall that Lebesgue almost everywhere
\begin{equation}\label{m function weight relation}
\frac{d\nu}{dx}(x) = \lim_{\eps \downarrow 0}\frac{1}{\pi}\im m(x
+ i\eps).
\end{equation}
Using this and the definition of $M$, it is easy to see that $M
\in \sm$ implies $d\nu \in \V$.

For the converse, assume that we can write
$$d\nu(x) = \sum_{j=1}^N c_j \delta(x-\lambda_j) + \bigl( \sqrt{2+x}
\bigr)^l \bigl( \sqrt{2-x} \bigr)^r v_0(x) dx$$ where $\lambda_j
\in \R\setminus\interval$, $c_j \in [0,1]$, $l,r \in \{\pm 1\}$,
and $\log v_0 \in \alg$.

Let $\{\psi_1,\psi_2\}$ be a partition of unity of $\interval$
subordinate to the cover $\{ [-2,1/2),(-1/2,2] \}$, with $\psi_1$
equalling one near $E=-2$ and $\psi_2$ equalling one near $E=2$.
Extend $\psi_1$ and $\psi_2$ to be zero outside $\interval$. In
this way we may write the $m$-function as
\begin{align*}
\begin{split}
m(z) &= \sum_{j=1}^N \frac{c_j}{\lambda_j -z} + \int
\frac{1}{x-z}\psi_1(x)d\nu(x) + \int
\frac{1}{x-z}\psi_2(x)d\nu(x)\\
&=s(z) + l(z) + r(z).
\end{split}
\end{align*}
By our choice of $\psi_1$ and $\psi_2$, and because $\lambda_j \in
\R \setminus \interval$, we have that $l(z)$ is smooth on $(1,2]$,
$r(z)$ is smooth on $[-2,-1)$, and $s(z)$ is smooth on
$\interval$.

We can now write $$M(z) = S(z) + L(z) + R(z)$$ where $$S(z) =
-s(z+z^{-1})$$ and similarly for $L$ and $R$.  Finally, we let
$$N(z) = M(z) - S(z) = L(z) + R(z).$$  As we have removed all the poles from
$M$, we see that $N$ is analytic in $\D$.  Moreover, because $S$
is smooth on $\partial \D$, it is clear that $M \in \sm$ if $N \in
\sm$, which we now prove.

We will first show that condition (1) holds for $N$.  Let $I_1$ be
an interval in $\partial \D$ avoiding $z = \pm 1$, and let $I$ a
slightly larger interval still avoiding $\pm 1$.  Let $\chi$ be a
smooth bump supported on $I$ equalling one on $I_1$.

By \eqref{m function weight relation} and the assumption that
$d\nu \in \V$ we see three things: $\sin (\theta) \im N(\theta)$
is a measure on $\partial \D$, $\chi(\theta)\sin(\theta)\im
N(\theta) \in \alg$, and $\im N$ is nonzero on $I$.

By Lemma \ref{lemma hilbert transform} we see
$$\chi(\theta) H[\sin(\theta) \im N(\theta)] = H[\chi(\theta)\sin(\theta)\im
N(\theta)] + f$$ where $f \in C^\infty$.  As
$\chi(\theta)\sin(\theta)\im N(\theta) \in \alg$ and $H$ is a
contraction in $\alg$ ($H$ multiplies the Fourier coefficients by
$0$ or $\pm i$, so it is a contraction in any space determined
only by Fourier coefficients), we see that
$$\chi(\theta) H[\sin(\theta) \im N(\theta)] \in \alg$$ too.  But
it is easy to see $(z - z^{-1})N(z)$ is analytic, so
$$H[\sin(\theta) \im N(\theta)] = \im (z-z^{-1})N(z) =
-\sin(\theta) \re N(\theta).$$ Combining these we find
$$\chi(\theta)\re N(\theta) = \chi(\theta) \frac{g(\theta)}{\sin
(\theta)}$$ for some $g \in \alg$.

Now we prove that (2) holds.  We will only consider the case $z_0
= 1$, the other case being similar.  Let $I_1$ be a $\partial
\D$-interval around $z_0$ to be chosen momentarily, let $I$ be a
slightly larger interval, and let $\chi$ be a smooth bump
supported on $I$ equalling one on $I_1$.

By \eqref{m function weight relation} and $d\nu \in \V$, we have
two cases to consider. Suppose first that
$$\chi(\theta) \im N(\theta) =
\chi(\theta)\frac{g(\theta)}{\sin(\theta)}$$ for some $g \in
\alg$.  Then arguing as in the proof of (1) shows $$\chi(\theta)
\re N(\theta) = \chi(\theta) \frac{h(\theta)}{\sin(\theta)}$$ for
some $h \in \alg$, so (2) holds in this case.

For the second case suppose $$\chi(\theta) \im N(\theta) =
\chi(\theta)\sin(\theta)g(\theta)$$ for some $g \in \alg$.  As $L$
is smooth near $z_0$ we have
\begin{align}\label{L form}
\begin{split}
\chi(\theta)\re L(\theta) &= \chi(\theta) \Biggl( \re L(0) +
\sin(\theta) \chi(\theta)\frac{\re L(\theta)-\re L(0)}{\sin (\theta)}\Biggr)\\
&=\chi(\theta)\Bigl( \re L(0) + \sin(\theta)h_1(\theta)\Bigr).
\end{split}
\end{align}
where $h_1 \in \alg$ if $I$ is chosen small enough.

Now we consider $\re R$ on $I$.  By assumption, $\im R$ is
continuous on $\partial \D$ and hence defines a measure. Also, $R$
is analytic in $\D$, so
\begin{align*}
\begin{split}
\chi(\theta) \re R(\theta) &= -\chi(\theta) H[\im R(\theta)]\\
&=H[-\chi(\theta)\im R(\theta)] + f_1\\
&=H[-\chi(\theta)\sin(\theta)g(\theta)] + H[\chi(\theta)\im L(\theta)] + f_1\\
&=\chi(\theta)\sin(\theta)H[-g] + H[\chi(\theta)\im L(\theta)]+ f_1  + f_2\\
&=\chi(\theta)\sin(\theta)h_2(\theta) + H[\chi(\theta)\im
L(\theta)]+ f_1 + f_2
\end{split}
\end{align*}
where the first equality follows from analyticity, the second and
fourth from Lemma \ref{lemma hilbert transform} (so $f_1, f_2 \in
C^\infty$), and the third from writing $R = M-L$.  As before, $h_2
\in \alg$ because $g \in \alg$.
 Because $L$ is smooth near $\theta = 0$, $\chi(\theta)\im
L(\theta)$ is smooth on all of $\partial \D$ if $I$ is chosen
small enough. In particular,
$$f := f_1 + f_2 + H[\chi(\theta)\im L(\theta)]$$ is smooth as
well.  Thus
\begin{align*}
\chi(\theta)f(\theta) &= \chi(\theta)\Biggl( f(0) +
\sin(\theta)\chi(\theta)\frac{f(\theta)-f(0)}{\sin(\theta)}\Biggr)\\
&=\chi(\theta)\Bigl( f(0) + \sin(\theta) h_3(\theta) \Bigr)
\end{align*}
with $h_3 \in \alg$, and so
\begin{equation}\label{R form}
\chi(\theta)\re R(\theta) = \chi(\theta)\Bigl( f(0) +
\sin(\theta)\bigl(h_2(\theta) + h_3(\theta)\bigr)\Bigr).
\end{equation}
Combining \eqref{L form} and \eqref{R form} shows (2) holds.
\end{proof}

%%%%%%%%%%%%%%%%%%%%%%%%%%%%%%%%%%%%%%%%%%%%%%%%%%%%%%%%%%%%%%%%%%%%%%%%%%%%%%%%%%%%%%%%%%%
%
%
%                                   Section
%
%
%%%%%%%%%%%%%%%%%%%%%%%%%%%%%%%%%%%%%%%%%%%%%%%%%%%%%%%%%%%%%%%%%%%%%%%%%%%%%%%%%%%%%%%%%%%

\section{$m$-functions and eigenvalues}
In this section we derive some properties of $\sm$.

\begin{proposition}\label{lemma m in M iff m1 in M}
Let $J$ a Jacobi matrix and let $J^{(1)}$ be the operator obtained
by removing the first row and column (from the top and left). Let
$M$ and $M^{(1)}$ be the M-functions corresponding to $J$ and
$J^{(1)}$.  Then $M \in \sm$ if and only if $M^{(1)} \in \sm$.
\end{proposition}

\begin{proof}
We will show that $M^{(1)} \in \sm$ implies $M \in \sm$, the other
direction being similar.

By \eqref{continued fraction for m} we have
\begin{gather*}
M(\theta) = \frac{1}{2 \cos \theta - b_1 - a_1^2
M^{(1)}(\theta)}\\
\im M = \Biggl| \frac{a_1}{2 \cos \theta - b_1 - a_1^2 M^{(1)}}
\Biggr|^2 \im M^{(1)}.
\end{gather*}
Let $I$ be an arc of $\partial \D$ missing $\theta = 0, \pi$.  As
$\im M^{(1)} \neq 0$ on $I$, we see
$$2 \cos \theta - b_1 - a_1^2 M^{(1)} \neq 0$$ on $I$.  By
assumption
$$2 \cos \theta - b_1 - a_1^2 M^{(1)} \in \algloc(I),$$ so by Theorem
\ref{lemma wiener levy} we have $\im M \in
\algloc(I)$ and is nonzero there, so part (1) of the definition of
$\sm$ holds.

Next we will show that if $M^{(1)}$ has the form required in part
(2) of the definition, then so does $M$. By hypothesis, we may
assume that $M^{(1)}(\theta) = c + (\sin \theta)^k g(\theta)$ on
some neighborhood $I$ of $\theta_0$, where $c \in \R$, $k \in
\{\pm 1\}$, and $g \in \algloc(I)$ with $\im g \neq 0$ there.

\underline{Case 1}: Suppose $k = -1$.  Then by subsuming the $c$
into $g$ we can write
\begin{align*}
M(\theta) &= \frac{1}{2 \cos \theta - b_1 - a_1^2 \frac{1}{\sin
\theta} g(\theta)}\\
&= \frac{\sin \theta}{(2 \cos \theta - b_1)(\sin \theta) - a_1^2
g(z)}\\
&=(\sin \theta) G(\theta).
\end{align*}
As $\alg$ is an algebra, the denominator of $G$ is in
$\algloc(I)$. As $\im g \neq 0$ we see that the denominator is
nonzero too.  By Theorem \ref{lemma wiener levy} we have $G \in
\algloc(I)$. As
$$\im G = \Biggl|\frac{1}{(2 \cos \theta - b_1)(\sin \theta) -
a_1^2 g(z)}\Biggr|^2 \im (-a_1^2 g)$$ we see $\im G \neq 0$ on
$I$.

\underline{Case 2}: Suppose $k = 1$.  Then we can write
\begin{align*}
M(\theta) &= \frac{1}{(2 \cos \theta - b_1 - a_1^2 c) - (\sin
\theta)a_1^2 g(\theta)}\\
&=\frac{1}{H(\theta) - (\sin \theta)a_1^2 g(\theta)}.
\end{align*}
If $H(\theta_0) = 0$, then because it is a real trigonometric
polynomial we can factor $H(\theta) = (\sin \theta)h(\theta)$ for
some real $h \in \algloc(I)$.  Then $h - a_1^2 g \in \algloc(I)$
and $\im (h - a_1^2 g) = \im (-a_1^2) \neq 0$, so
$$M(\theta) = \frac{1}{\sin \theta} \frac{1}{h - a_1^2 g} =
\frac{1}{\sin \theta} G(\theta)$$ where $G \in \algloc(I)$ and has
nonvanishing imaginary part.

If $H(\theta_0) = \tfrac{1}{c}$ for some constant $c \in \R
\setminus \{0\}$, then as before we can write $H(\theta_0) -
H(\theta) = (\sin \theta)h(\theta)$ for some real $h \in
\algloc(I)$.  Then
\begin{align*}
M(\theta) - c &= \frac{1}{H(\theta) - (\sin \theta) a_1^2
g(\theta)}-\frac{1}{H(\theta_0)}\\
&= c \frac{\bigl( H(\theta_0) - H(\theta)\bigr) - (\sin \theta)
a_1^2 g(\theta)}{H(\theta) - (\sin \theta) a_1^2 g(\theta)}\\
&= (\sin \theta) c \frac{h(\theta) - a_1^2 g(\theta)}{H(\theta) -
(\sin \theta) a_1^2 g(\theta)}\\
&= (\sin \theta) G(\theta).
\end{align*}
Both the numerator and the denominator of $G$ are in $\algloc(I)$.
The denominator is nonvanishing in a neighborhood $I' \subseteq I$
of $\theta_0$. Thus, $G \in \algloc (I')$ by Theorem \ref{lemma
wiener levy}. Next we compute
\begin{equation*}
\im G (\theta) =c \frac{-a_1^2 \im g \Bigl( \re H - (\sin \theta)
a_1^2 \re g \Bigr)
 + a_1^2 (\sin \theta) \im g \Bigl( \re h - a_1^2
\re g \Bigr)}{|H - (\sin \theta)a_1^2 g|^2}.
\end{equation*}
The denominator is nonvanishing and in $\algloc(I')$.  The second
term in the numerator vanishes at $\theta_0$, but the first tends
to $-a_1^2 \im g(\theta_0) H(\theta_0) \neq 0$.  So $\im G \neq 0$
in some neighborhood $I'' \subseteq I'$ of $\theta_0$.
\end{proof}

\begin{proposition}\label{lemma m in M iff tilde m in M}
Suppose $J$ and $\widetilde J$ are two Jacobi matrices satisfying
$J^{(1)} = \widetilde J^{(1)}$.  Then $M \in \sm$ if and only if
$\widetilde M \in \sm$.
\end{proposition}

\begin{proof}
This follows immediately from two applications of Proposition
\ref{lemma m in M iff m1 in M}.
\end{proof}

\begin{proposition}\label{lemma m no eits implies tilde m no eits}
Let $J$ be a Jacobi matrix with no eigenvalues off $[-2,2]$, and
assume that $M \in \sm$. Then there is a unique doubly-resonant
$\widetilde J$ with $\widetilde J ^{(1)} = J^{(1)}$ (so
$\widetilde M \in \sm$) and no eigenvalues off $[-2,2]$.
\end{proposition}

\begin{proof}
By \eqref{continued fraction for m} we have $$M(\theta) =
\frac{1}{2 \cos \theta - b_1 - a_1^2 M^{(1)}(\theta)}.$$
Similarly, if $\widetilde J ^{(1)} = J^{(1)}$ then we have
$$\widetilde M (\theta) = \frac{1}{2 \cos \theta -\tilde b_1 - \tilde
a_1^2 \widetilde M^{(1)}(\theta)} = \frac{1}{2 \cos \theta -
\tilde b_1 - \tilde a_1^2 M^{(1)}(\theta)}.$$ Combining these one
finds
\begin{equation}\label{m tilde m relationship}
\widetilde M(\theta) = \frac{a_1^2}{\tilde a_1^2 M (\theta)^{-1} -
\Delta (\theta)}
\end{equation}
where $\Delta (\theta) = (\delta a) (2 \cos \theta) - \delta ab$,
$\delta a = \tilde a_1^2 - a_1^2$, and $\delta ab = \tilde a_1^2
b_1 - a_1^2 \tilde b_1$.

As $M \in \sm$, in some $\partial \D$-neighborhood $I_+$ of
$\theta = 0$ we can write $$M(\theta) = c_+ + (\sin
\theta)^{k_+}g_+(\theta)$$ for some $c_+ \in \R$, $k_+ \in \{\pm
1\}$, and $g_+ \in \algloc(I_+)$ with $\im g_+ \neq 0$. Similarly,
in a neighborhood $I_-$ of $\theta = \pi$ we can write $$M(\theta)
= c_- + (\sin \theta)^{k_-}g_-(\theta).$$

By \eqref{m tilde m relationship} we see that to make $\widetilde
J$ doubly-resonant, we must choose $\tilde a_1$ and $\tilde b_1$
so that $\tilde a_1^2 M(\theta)^{-1} - \Delta(\theta) = 0$ at
$\theta = 0, \pi$. There are four cases depending on the various
combinations of $k_-$ and $k_+$. When $k_- = k_+ = -1$ we just
choose $\tilde a_1 = a_1$ and $\tilde b_1 = b_1$.  When $k_- = 1$
and $k_+ = -1$ choose
$$\tilde a_1^2 = a_1^2 \Bigl(\frac{4c_-}{4c_- + 1}\Bigr) \quad \text{and}\quad
\tilde b_1 = 2 \Bigl(\frac{2b_1 c_- + 1}{4c_- + 1}\Bigr).$$ When
$k_- = -1$ and $k_+ = 1$ choose
$$\tilde a_1^2 = a_1^2 \Bigl(\frac{4c_+}{4c_+ - 1}\Bigr) \quad \text{and}\quad
\tilde b_1 = 2 \Bigl(\frac{2b_1 c_+ + 1}{4c_+ - 1}\Bigr).$$ When
$k_- = k_+ = 1$ choose
$$\tilde a_1^2 = a_1^2 \Bigl(\frac{4c_-c_+}{4c_-c_+ - c_- + c_+}\Bigr) \quad \text{and}\quad
\tilde b_1 = 2 \Bigl(\frac{2b_1 c_-c_+ + c_- + c_+}{4c_-c_+ - c_-
+ c_+}\Bigr).$$

Of course, we must check that $\tilde a_1^2 > 0$ so that
$\widetilde J$ really is a Jacobi matrix.  This amounts to showing
that $c_- < -1/4$ and $c_+ > 1/4$.  As $J$ has no eigenvalues off
$[-2,2]$, $$m(E) = \int_{-2}^2 \frac{d\nu(x)}{x-E}$$ where $d\nu$
is the spectral measure corresponding to $J$.  For $t \in [-2,2]$
and $E > 2$, $t-E \geq -4$. So because $d\nu$ is a probability
measure that is not a point mass at $t=2$ we have
$$M(\theta = 0) = \lim_{E \downarrow 2}-m(E) > 1/4.$$ Similar
arguments show $M(\theta = \pi) < -1/4$.

It remains to show that $\widetilde J$ has no eigenvalues off
$[-2,2]$, or equivalently that $\widetilde M$ has no poles on
$(-1,1)$.  As $J$ has no eigenvalues off $[-2,2]$, $M$ is analytic
on $\D$.  So by \eqref{m tilde m relationship} it suffices to show
that $f(E) := \tilde a_1^2 + m(E)\bigl(\delta a E - \delta ab
\bigr) \neq 0$ for $|E| > 2$.  As $J$ has no eigenvalues off
$[-2,2]$ we have
$$\frac{d}{dE}m(E) = \int_{-2}^2 \frac{d\nu(x)}{(x-E)^2} > 0.$$
Since $\delta a E - \delta ab$ is linear in $E$, $f$ is monotone
in $E$ for $|E| > 2$.  By our choice of $\tilde a_1^2$ and $\tilde
b_1$ we have $f(\pm 2) = 0$, and so $f(E) \neq 0$ for $|E| > 2$,
as required.
\end{proof}

%%%%%%%%%%%%%%%%%%%%%%%%%%%%%%%%%%%%%%%%%%%%%%%%%%%%%%%%%%%%%%%%%%%%%%%%%%%%%%%%%%%%%%%%%%%
%
%
%                                   Section
%
%
%%%%%%%%%%%%%%%%%%%%%%%%%%%%%%%%%%%%%%%%%%%%%%%%%%%%%%%%%%%%%%%%%%%%%%%%%%%%%%%%%%%%%%%%%%%

\section{Proof of Theorem \ref{lemma main thm alg case}}
As a final preliminary, we recall the definition of the
Carath\'eodory function associated to $d\mu = w
\frac{d\theta}{2\pi}$:
$$F(z) = \int_0^{2\pi}\frac{e^{i\theta}+z}{e^{i\theta}-z}d\mu
(\theta).$$ Note that almost everywhere $$\lim_{r \uparrow 1}\re
F(r e^{i \theta}) = w(\theta)$$ and if $M$ has no poles in $\D$
then
\begin{equation}\label{M F relation}
M(z) = \frac{-F(z)}{z - z^{-1}}.
\end{equation}

\begin{proof}[Proof of Theorem \ref{lemma main thm alg case}]
Suppose $\lambda, \kappa \in \seqalg$.  Choose $N$ so large that
$\lambda(J^{(N)}), \kappa(J^{(N)}) \in \seqalg$ with small enough
norms to apply Proposition \ref{lemma solves Geronimus}. This
produces a sequence $\tilde \alpha \in \seqalg$ solving
\eqref{expanded lambda kappa}, but may not have $\tilde
\alpha_{-1} = -1$.  If we change $\tilde \alpha_{-1}$ to be $-1$,
this changes the top row and column of $J^{(N)}$ producing a
matrix $\widetilde J$. As $\widetilde J$ has a sequence of
Verblunsky parameters (namely $\tilde \alpha$), we see $\sigma
(\widetilde J) \subseteq \interval$.

Let $\widetilde M$ be the $M$-function associated to $\widetilde
J$. We will show $\widetilde M \in \sm$. Let $d\tilde\mu$ be the
measure on $\partial \D$ corresponding to $\tilde \alpha$.  As
$\tilde \alpha \in \seqalg$, we may apply Theorems \ref{lemma
Golinskii-Ibragimov} and \ref{lemma baxter} to find $d\tilde\mu =
\tilde w \frac{d\theta}{2\pi}$ and $\log \tilde w \in \alg$.  As
$\tilde w = \re \widetilde F$ we see that $\log (\re F) \in \alg$
too, so by Theorem \ref{lemma wiener levy}, $\re F \in \alg$ and
is nonvanishing. As $\re f \mapsto \im f$ is a contraction in
$\alg$ we see $\widetilde F \in \alg$ and is nonvanishing.  By
\eqref{M F relation} we have
\begin{gather*}
\widetilde M(\theta) = i \frac{1}{\sin \theta} \frac{\widetilde
F(\theta)}{2}\\
\im \widetilde M = \frac{1}{\sin \theta} \frac{\tilde
w(\theta)}{2}.
\end{gather*}
In particular, $\widetilde M \in \sm$, as claimed.

By Proposition \ref{lemma m in M iff tilde m in M}, $M(J^{(N)})
\in \sm$, and by repeated applications of Proposition \ref{lemma m
in M iff m1 in M} we have that $M \in \sm$. Finally, $J$ and
$\widetilde J$ differ by a finite-rank perturbation.  Since
$\widetilde J$ has no eigenvalues off $[-2,2]$ and a finite-rank
perturbation can only produce a finite number of eigenvalues in
each spectral gap, $J$ has only finitely-many eigenvalues and they
all lie in $\R \setminus [-2,2]$. By Proposition \ref{lemma dnu M
connection} we have $d\nu \in \V$.

Now consider the converse.  As $J$ has only finitely-many
eigenvalues off $[-2,2]$, the Sturm Oscillation Theorem guarantees
we can choose $N$ large enough that $J^{(N)}$ has no eigenvalues
off $[-2,2]$. By Propositions \ref{lemma m in M iff m1 in M},
\ref{lemma m in M iff tilde m in M}, and \ref{lemma m no eits
implies tilde m no eits}, there is a unique doubly-resonant Jacobi
matrix $\widetilde J$ with $\widetilde J^{(1)} = J^{(N+1)}$,
$\widetilde M \in \sm$, and no eigenvalues off $[-2,2]$.  As
above,
$$\widetilde M = i\frac{1}{\sin \theta}\frac{\widetilde
F(\theta)}{2}$$ so $\widetilde F \in \alg$ and $\widetilde w$ is
nonvanishing. By Theorems \ref{lemma Golinskii-Ibragimov} and
\ref{lemma baxter} we have $\widetilde \alpha \in \seqalg$, where
$\widetilde \alpha$ is the sequence of Verblunsky parameters
corresponding to $\widetilde J$. But then $\lambda(\widetilde
J),\kappa(\widetilde J) \in \seqalg$ by Lemma \ref{lemma alpha
implies parameters}.  As $\lambda(J)$ and $\kappa(J)$ differ from
$\lambda(\widetilde J)$ and $\kappa(\widetilde J)$ by only
finitely-many terms, we have $\lambda(J),\kappa(J) \in \seqalg$
too.
\end{proof}

%%%%%%%%%%%%%%%%%%%%%%%%%%%%%%%%%%%%%%%%%%%%%%%%%%%%%%%%%%%%%%%%%%%%%%%%%%%%%%%%%%%%%%%%%%%
%
%
%                                 Bibliography
%
%
%%%%%%%%%%%%%%%%%%%%%%%%%%%%%%%%%%%%%%%%%%%%%%%%%%%%%%%%%%%%%%%%%%%%%%%%%%%%%%%%%%%%%%%%%%%

\end{document}